\UseRawInputEncoding
\documentclass[oneside, 12pt]{amsart}
\usepackage{amsmath}
\usepackage{amssymb}
\usepackage{dsfont}
\usepackage{color}
\usepackage[usenames,dvipsnames,x11names,svgnames]{xcolor}
\usepackage[colorlinks=true,linkcolor=NavyBlue,citecolor=DarkGreen, urlcolor=Blue]{hyperref}
\usepackage{amsthm}
\usepackage{amscd}
\usepackage{geometry}
\usepackage{mathrsfs}
\usepackage{mathtools}
\usepackage{esint}

\theoremstyle{definition}

\newcommand{\be}{\begin{equation*}}
\newcommand{\ee}{\end{equation*} }

\newcommand{\ben}{\begin{equation}}
\newcommand{\een}{\end{equation} }

\newcommand{\bs}{\begin{split}}
\newcommand{\es}{\end{split}}

\newcommand{\bmu}{\begin{multline*}}
\newcommand{\emu}{\end{multline*}}

\newcommand{\bmun}{\begin{multline}}
\newcommand{\emun}{\end{multline}}

\begin{document}



\title[]{A scaling law chaotic system}

\author{Xiao-Jun Yang}

\email{dyangxiaojun@163.com; xjyang@cumt.edu.cn}

\address{School of Mathematics, China University of Mining and Technology, Xuzhou 221116, China}

\begin{abstract}
In this article, we propose an anomalous chaotic system of the scaling-law ordinary
differential equations involving the Mandelbrot scaling law. This chaotic behavior
shows the "Wukong" effect. The comparison among the Lorenz and scaling-law attractors is
discussed in detail. We also suggest the conjecture for the fixed point theory for the fractal SL attractor.
The scaling-law chaos may be open a new door in the study of the chaos theory.
\end{abstract}

\maketitle

\section{Introduction}\label{sec:1}
The Lorenz attractor, discovered in 1963 by the MIT meteorologist Edward Lorenz to describe a
simplified mathematical model for atmospheric convection, is a chaotic
system of three ordinary differential equations \cite{1}:
\begin{equation}
\label{eq1}
\frac{dx\left( t \right)}{dt}=10\left( {y\left( t \right)-x\left( t \right)}
\right),
\end{equation}
\begin{equation}
\label{eq2}
\frac{dy\left( t \right)}{dt}=x\left( t \right)\left( {\frac{8}{3}-z\left( t
\right)} \right)-y\left( t \right),
\end{equation}
and
\begin{equation}
\label{eq3}
\frac{dz\left( t \right)}{dt}=x\left( t \right)y\left( t \right)-28z\left( t
\right),
\end{equation}
where $10$ is the Prandtl number, $8/3$ is proportional to the~Rayleigh
number, and $28$ is a geometric factor. The Lorenz equations (\ref{eq1}), (\ref{eq2}) and
(\ref{eq3}) are used to study the mathematical models in lasers \cite{2}, dynamos \cite{3},
thermosyphons \cite{4}, DC motors \cite{5}, electric circuits \cite{6}, chemical
reactions \cite{7}, and  waterwheel \cite{8,9}.

In the present article, we consider that the fractal scaling-law (SL)
chaotic system is given by the SL ordinary differential equations
\begin{equation}
\label{eq4}
{ }^{MSL}D_t^{\left( 1 \right)} x\left( t \right)=a\left( {y\left( t
\right)-x\left( t \right)} \right),
\end{equation}
\begin{equation}
\label{eq5}
{ }^{MSL}D_t^{\left( 1 \right)} y\left( t \right)=x\left( t \right)\left(
{b-z\left( t \right)} \right)-y\left( t \right),
\end{equation}
and
\begin{equation}
\label{eq6}
{ }^{MSL}D_t^{\left( 1 \right)} z\left( t \right)=x\left( t \right)y\left( t
\right)-cz\left( t \right),
\end{equation}
in which the fractal SL derivative of the function $\aleph \left( t \right)$
is defined as \cite{10,11,11a,11b}
\begin{equation}
\label{eq7}
{ }^{MSL}D_t^{\left( 1 \right)} \aleph \left( t \right)=\frac{t^D}{\lambda
}\frac{d\aleph \left( t \right)}{dt},
\end{equation}
where $a$, $b$, and $c$ are the given parameters, the coefficient $\lambda $
is given by
\begin{equation}
\label{eq8}
\lambda =\mu \left( {1-D} \right),
\end{equation}
the Mandelbrot scaling law $\Lambda \left( {\mu ,D,t} \right)$ is expressed
in the form \cite{12}
\begin{equation}
\label{eq9}
\Lambda \left( {\mu ,D,t} \right)=\mu t^{1-D}
\end{equation}
for the parameter $\mu >0$, time $t$, and fractal dimension $0<D<1$.

The main of the paper is to study the fractal SL chaotic system which are
given by the fractal SL ordinary differential equations, and to present the "Wukong effect", which is
observed in the plots of the plane x-y for the fractal SL attractors. The structure of the paper is designed as follows.
In Section \ref{Sec:2} we propose a fractal SL attractor with the variable parameter.
In Section \ref{Sec:3} we observe the typical systems for the fractal SL ordinary differential
equations. In Section \ref{Sec:4} we present the comparative results among the Lorenz and fractal SL attractors.
The conclusion and future Work are given in Section \ref{Sec:5}.

\section{The fractal SL attractor with the variable parameter} \label{Sec:2}

By using (\ref{eq4}), (\ref{eq5}), and (\ref{eq6}), the fractal SL attractor is represented by the fractal SL
ordinary differential equations:
\begin{equation}
\label{eq10}
\frac{10}{3}t^{\frac{2}{3}}\frac{dx\left( t \right)}{dt}=a\left( {y\left( t
\right)-x\left( t \right)} \right),
\end{equation}
\begin{equation}
\label{eq11}
\frac{10}{3}t^{\frac{2}{3}}\frac{dy\left( t \right)}{dt}=x\left( t
\right)\left( {\frac{3}{10}-z\left( t \right)} \right)-y\left( t \right),
\end{equation}
and
\begin{equation}
\label{eq12}
\frac{10}{3}t^{\frac{2}{3}}\frac{dz\left( t \right)}{dt}=x\left( t
\right)y\left( t \right)-27z\left( t \right),
\end{equation}
where $a$ is a real variable parameter.

\begin{figure}[h]
\centering
\includegraphics[width=9cm]{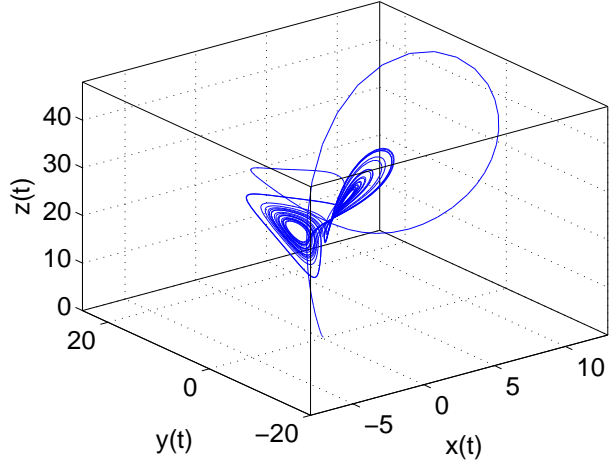}
\caption{The fractal SL ordinary differential equations with the parameter $a =2.35$.}
\label{Fig1}
\end{figure}

\begin{figure}[h]
\centering
\includegraphics[width=9cm]{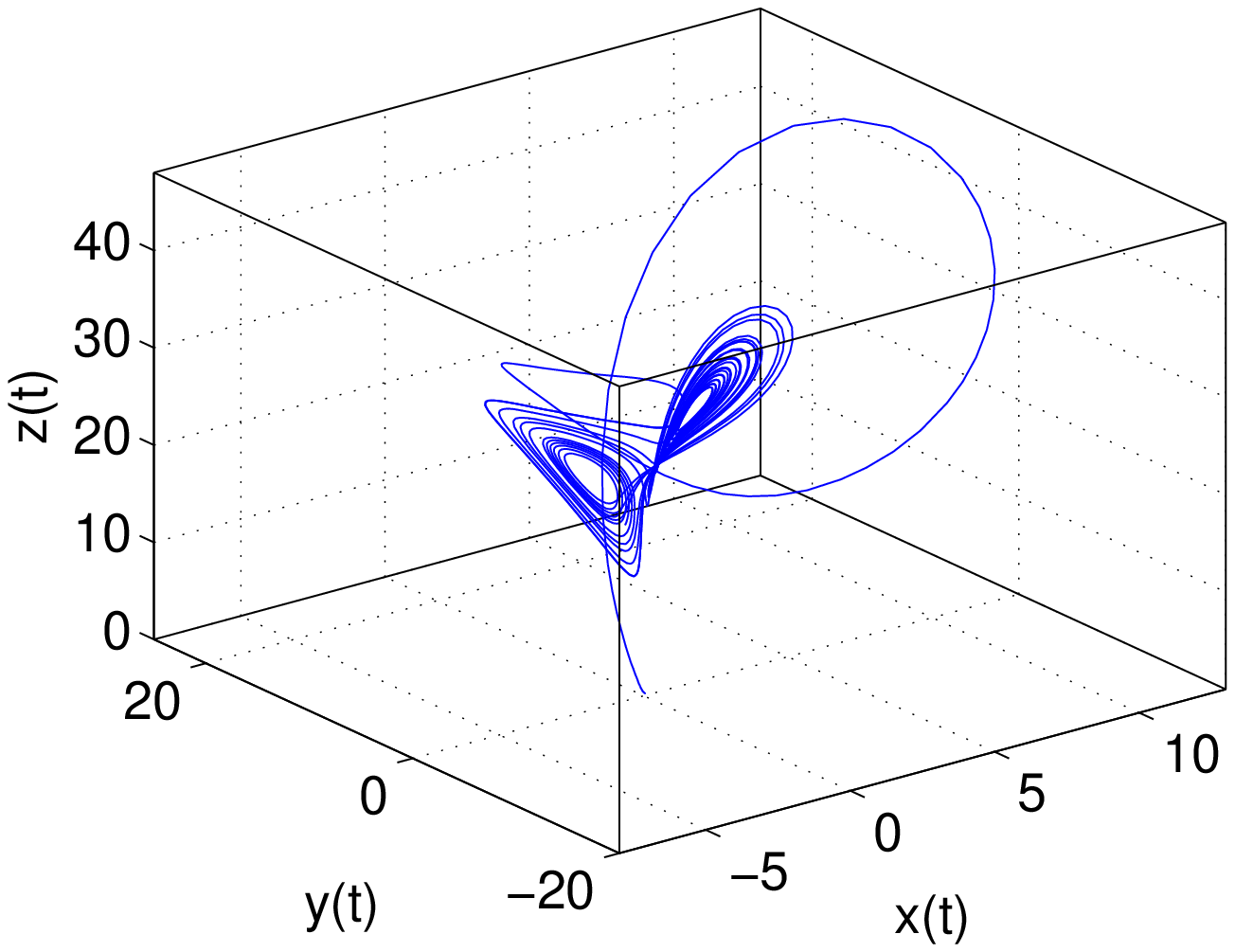}
\caption{The fractal SL ordinary differential equations with the parameter $a =2$.}
\label{Fig2}
\end{figure}

\begin{figure}[h]
\centering
\includegraphics[width=9cm]{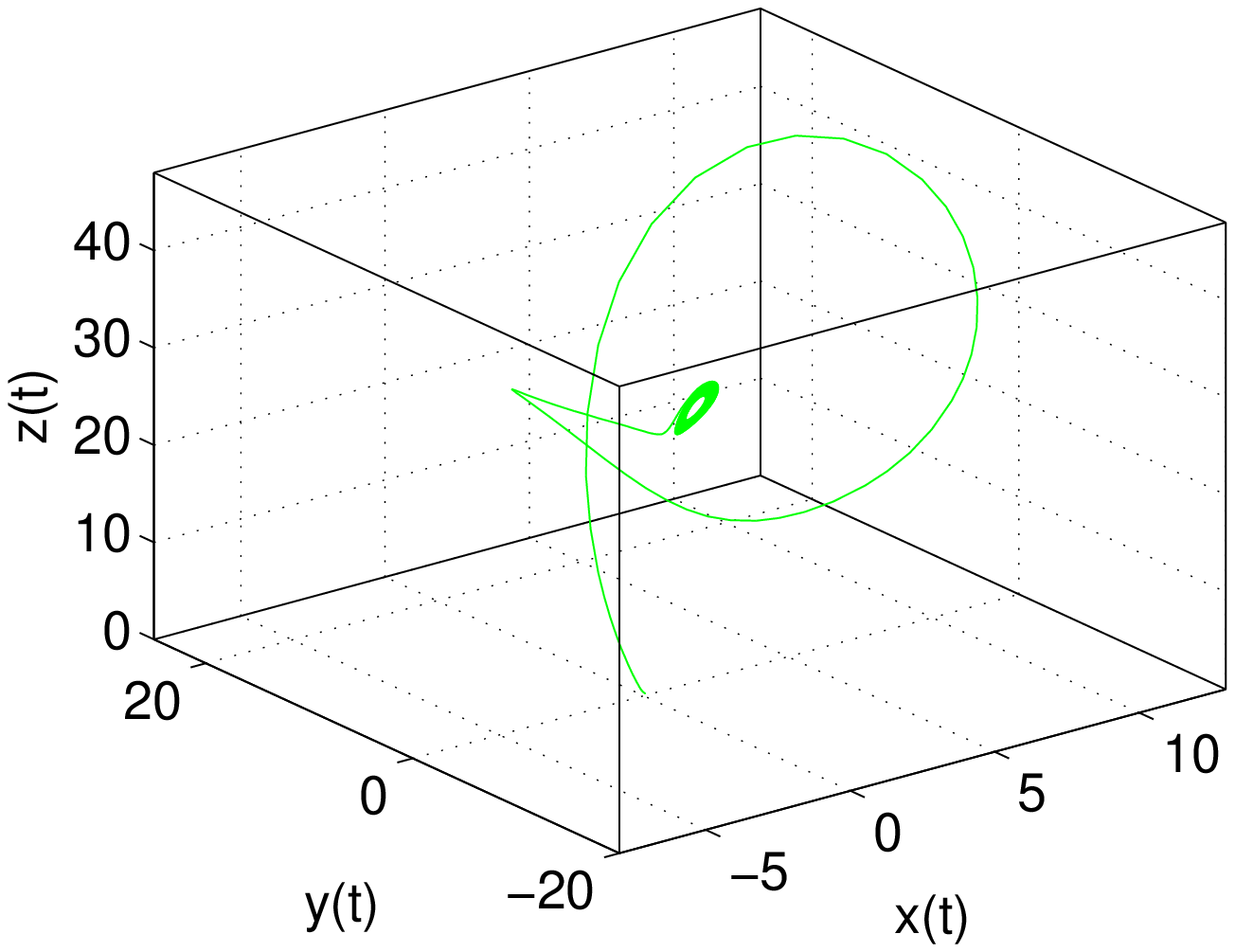}
\caption{The fractal SL ordinary differential equations with the parameter $a =1.5$.}
\label{Fig3}
\end{figure}

\begin{figure}[h]
\centering
\includegraphics[width=9cm]{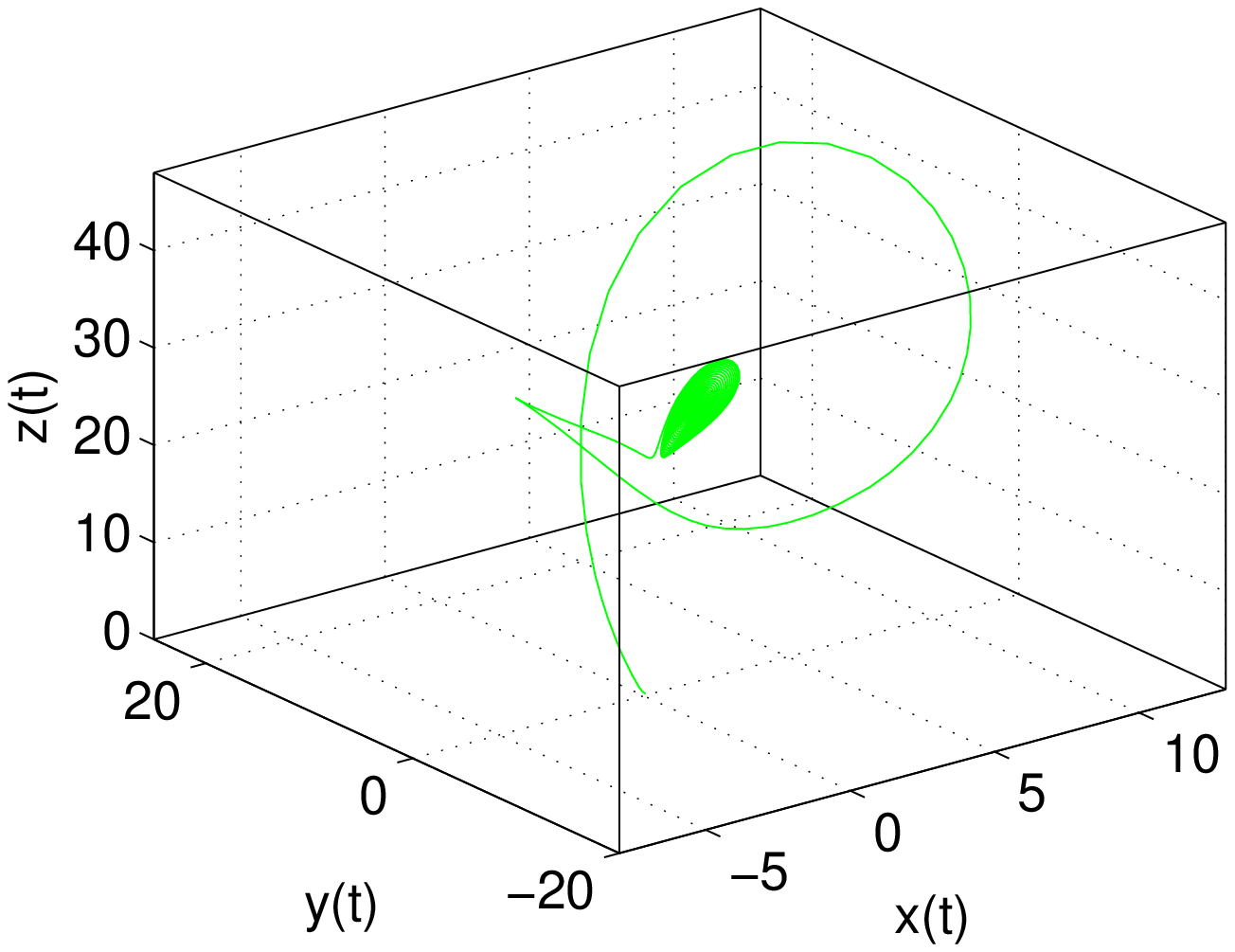}
\caption{The fractal SL ordinary differential equations with the parameter $a =1.35$.}
\label{Fig4}
\end{figure}

The anomalous behaviors of the fractal SL attractor with the variable parameter $a$
and the parameters $D=2/3$ and $\mu =0.9$ are considered when the
initial conditions are $x\left( 0 \right)=0.1$, $y\left( 0 \right)=0.1$, and
$z\left( 0 \right)=0.1$, and time $t$ changes from $0.1$ to $10^6$.
Figure \ref{Fig1} shows the fractal SL ordinary differential equations with the parameters
$a=2.35$, $D=2/3$ and $\mu =0.9$.
The fractal SL ordinary differential equations with the parameters
$a=2.35$, $D=2/3$ and $\mu =0.9$ is given in Figure \ref{Fig2}.
The fractal SL ordinary differential equations with the parameters
$a=1.5$, $D=2/3$ and $\mu =0.9$ is presented in Figure \ref{Fig3}.
The fractal SL ordinary differential equations with the parameters
$a=1.35$, $D=2/3$ and $\mu =0.9$ is illustrated in Figure \ref{Fig4}.
It is observed that the fractal SL ordinary differential equations with the parameters
$D=2/3$ and $\mu =0.9$ have the chaotic behaviors when $a=2.35$ and $a=2$.

\section{The fractal SL attractors} \label{Sec:3}

To present the anomalous behaviors of the fractal SL attractors, we start to
observe the typical systems for the fractal SL ordinary differential
equations.

\subsection{The fractal SL attractor I}

\begin{figure}[h]
\centering
\includegraphics[width=9cm]{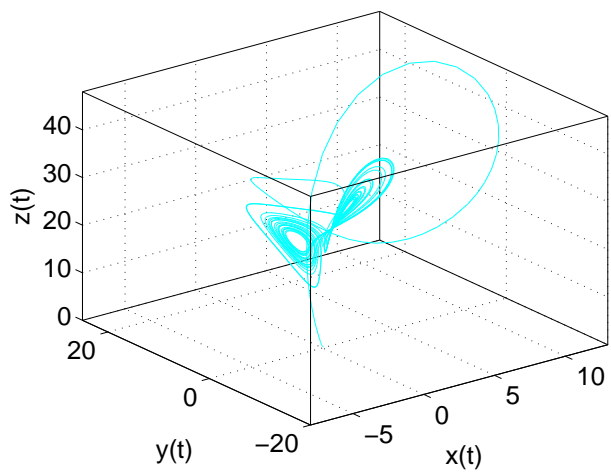}
\caption{The fractal SL attractor I with the parameters $a =2.35$, $D=2/3$
and $\mu =0.9$.}
\label{Fig5}
\end{figure}

Making use of (\ref{eq10}), (\ref{eq11}), and (\ref{eq12}) with the parameter $D=2/3$ and $\mu
=0.9$, the fractal SL attractor I for the fractal SL differential equations
can be presented as follows:
\begin{equation}
\label{eq13}
\frac{10}{3}t^{\frac{2}{3}}\frac{dx\left( t \right)}{dt}=\frac{47}{20}\left(
{y\left( t \right)-x\left( t \right)} \right),
\end{equation}
\begin{equation}
\label{eq14}
\frac{10}{3}t^{\frac{2}{3}}\frac{dy\left( t \right)}{dt}=x\left( t
\right)\left( {\frac{3}{10}-z\left( t \right)} \right)-y\left( t \right),
\end{equation}
and
\begin{equation}
\label{eq15}
\frac{10}{3}t^{\frac{2}{3}}\frac{dz\left( t \right)}{dt}=x\left( t
\right)y\left( t \right)-27z\left( t \right).
\end{equation}
The chaotic trajectories of the SL attractor I with the parameters $D=2/3$
and $\mu =0.9$ are showed in Figure \ref{Fig5}, where the initial conditions are
$x\left( 0 \right)=0.1$, $y\left( 0 \right)=0.1$, and $z\left( 0
\right)=0.1$, and time $t$ changes from $0.1$ to $10^6$.

\begin{figure}[h]
\centering
\includegraphics[width=9cm]{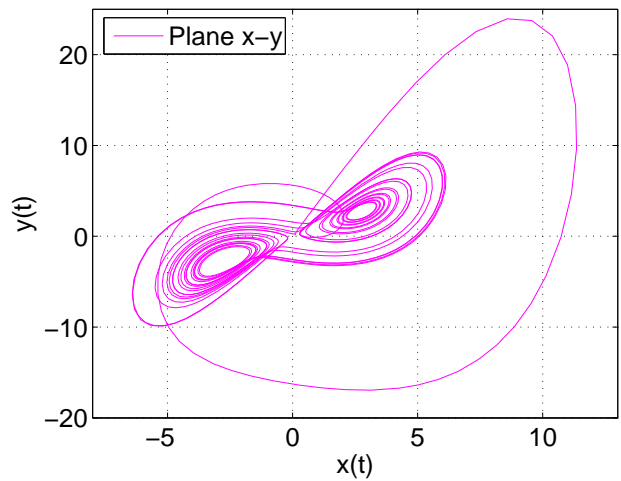}
\caption{The plane x-y with the parameters $a =2.35$, $D=2/3$
and $\mu =0.9$.}
\label{Fig6}
\end{figure}

The plane x-y of the fractal SL attractor I is shown in Figure \ref{Fig6},
where the initial conditions are $x\left( 0 \right)=0.1$, $y\left( 0
\right)=0.1$, and $z\left( 0 \right)=0.1$, and time $t$ changes from $0.1$
to $10^6$.

\begin{figure}[h]
\centering
\includegraphics[width=9cm]{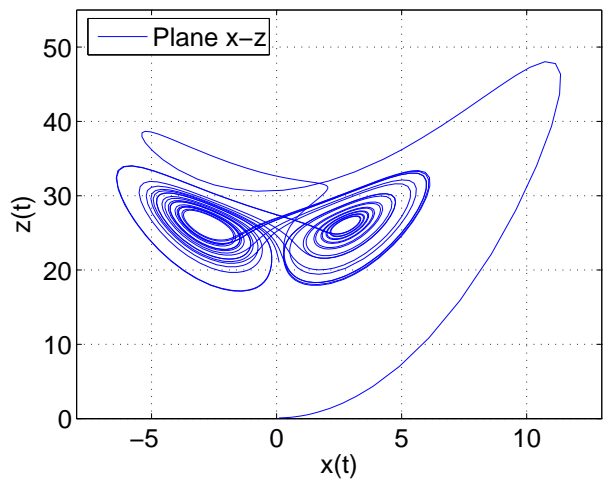}
\caption{The plane x-z with the parameters $a =2.35$, $D=2/3$
and $\mu =0.9$.}
\label{Fig7}
\end{figure}

The plane x-z of the fractal SL attractor I is given in Figure \ref{Fig7},
where the initial conditions are $x\left( 0 \right)=0.1$, $y\left( 0
\right)=0.1$, and $z\left( 0 \right)=0.1$, and time $t$ changes from $0.1$
to $10^6$.

\begin{figure}[h]
\centering
\includegraphics[width=9cm]{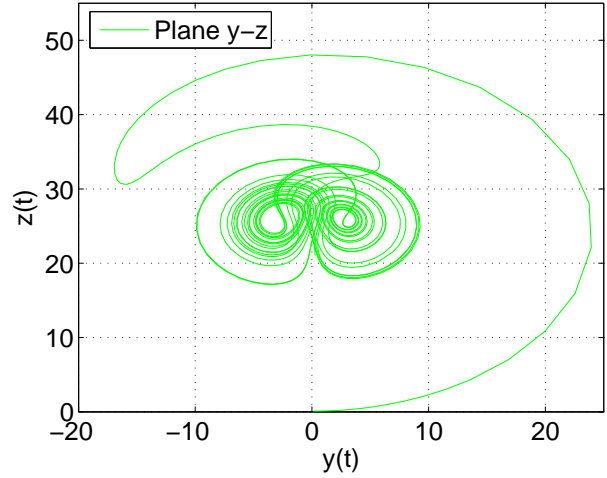}
\caption{The plane y-z with the parameters $a =2.35$, $D=2/3$
and $\mu =0.9$.}
\label{Fig8}
\end{figure}

The plane y-z of the fractal SL attractor I is depicted in Figure \ref{Fig8},
where the initial conditions are $x\left( 0 \right)=0.1$, $y\left( 0
\right)=0.1$, and $z\left( 0 \right)=0.1$, and time $t$ changes from $0.1$
to $10^6$.

\begin{figure}[h]
\centering
\includegraphics[width=10.5cm]{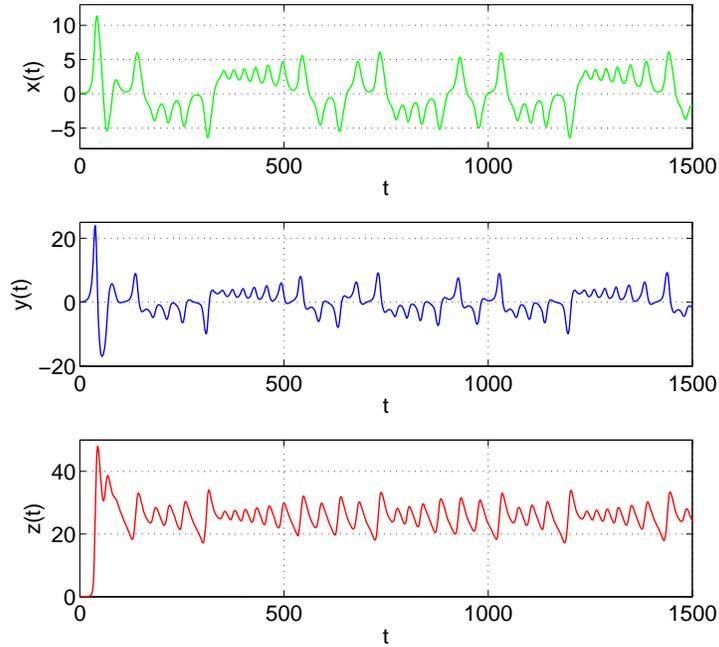}
\caption{The time series for $x\left( t \right)$, $y\left( t \right)$,
and $z\left( t \right)$ with the parameters $a =2.35$, $D=2/3$
and $\mu =0.9$.}
\label{Fig9}
\end{figure}

The plots of the time series for $x\left( t \right)$, $y\left( t \right)$,
and $z\left( t \right)$ are plotted in Figure \ref{Fig9}, where $D=2/3$, $\mu =0.9$,
and time $t$ changes from $0.1$ to $10^6$.

\subsection{The fractal SL attractor II}

Similarly, by using (\ref{eq10}), (\ref{eq11}), and (\ref{eq12}) with the parameter $D=2/3$ and
$\mu =0.9$, we obtain the
\begin{equation}
\label{eq16}
\frac{10}{3}t^{\frac{2}{3}}\frac{dx\left( t \right)}{dt}=2\left( {y\left( t
\right)-x\left( t \right)} \right),
\end{equation}
\begin{equation}
\label{eq17}
\frac{10}{3}t^{\frac{2}{3}}\frac{dy\left( t \right)}{dt}=x\left( t
\right)\left( {\frac{3}{10}-z\left( t \right)} \right)-y\left( t \right),
\end{equation}
and
\begin{equation}
\label{eq18}
\frac{10}{3}t^{\frac{2}{3}}\frac{dz\left( t \right)}{dt}=x\left( t
\right)y\left( t \right)-27z\left( t \right).
\end{equation}

\begin{figure}[h]
\centering
\includegraphics[width=9cm]{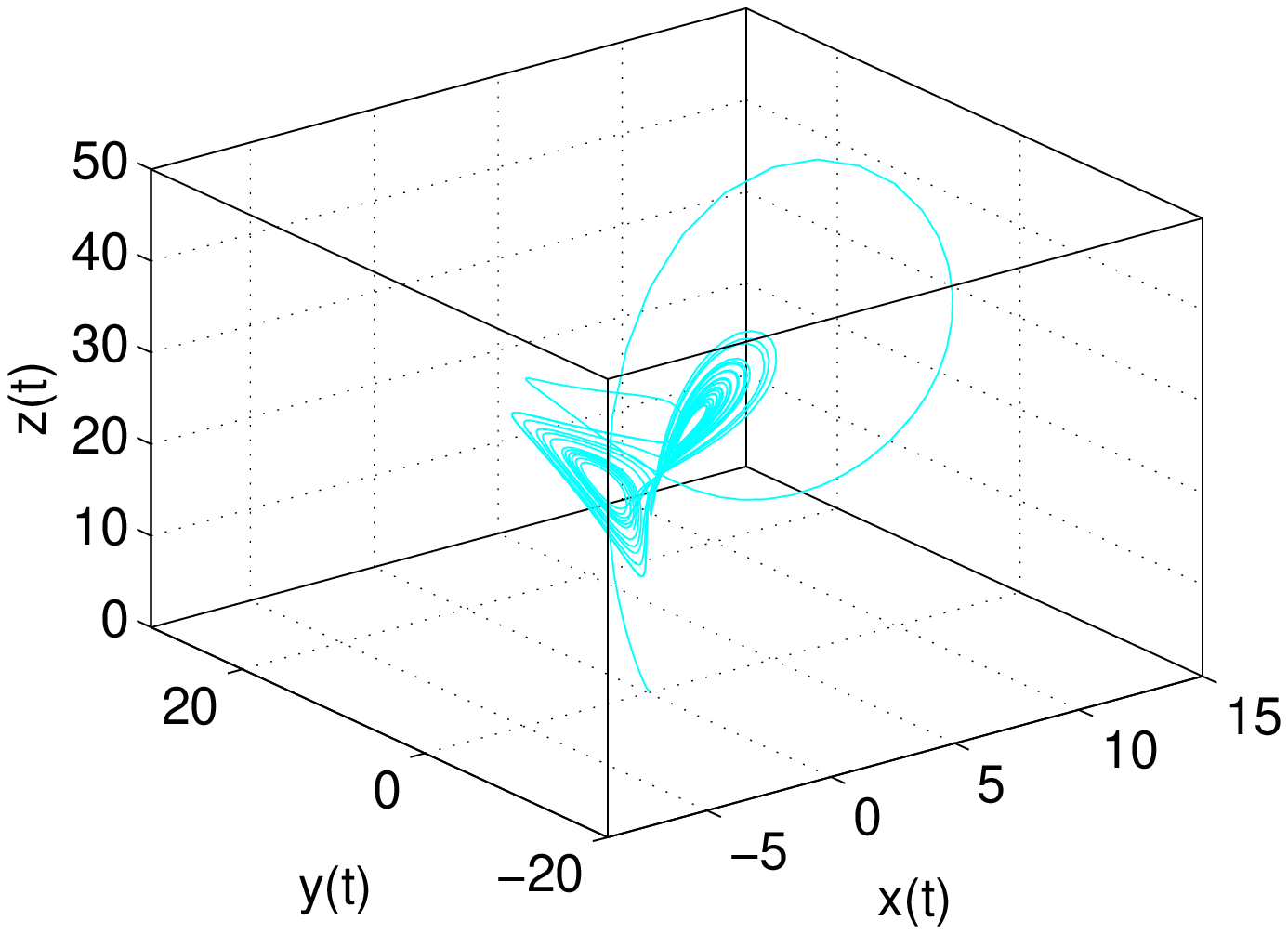}
\caption{The fractal SL attractor II with the parameters $a=2$, $D=2/3$
and $\mu =0.9$.}
\label{Fig10}
\end{figure}

The chaotic trajectories of the SL attractor II with the parameters $D=2/3$
and $\mu =0.9$ are showed in Figure \ref{Fig10}, where the initial conditions are
$x\left( 0 \right)=0.1$, $y\left( 0 \right)=0.1$, and $z\left( 0
\right)=0.1$, and time $t$ changes from $0.1$ to $10^6$.

\begin{figure}[h]
\centering
\includegraphics[width=9cm]{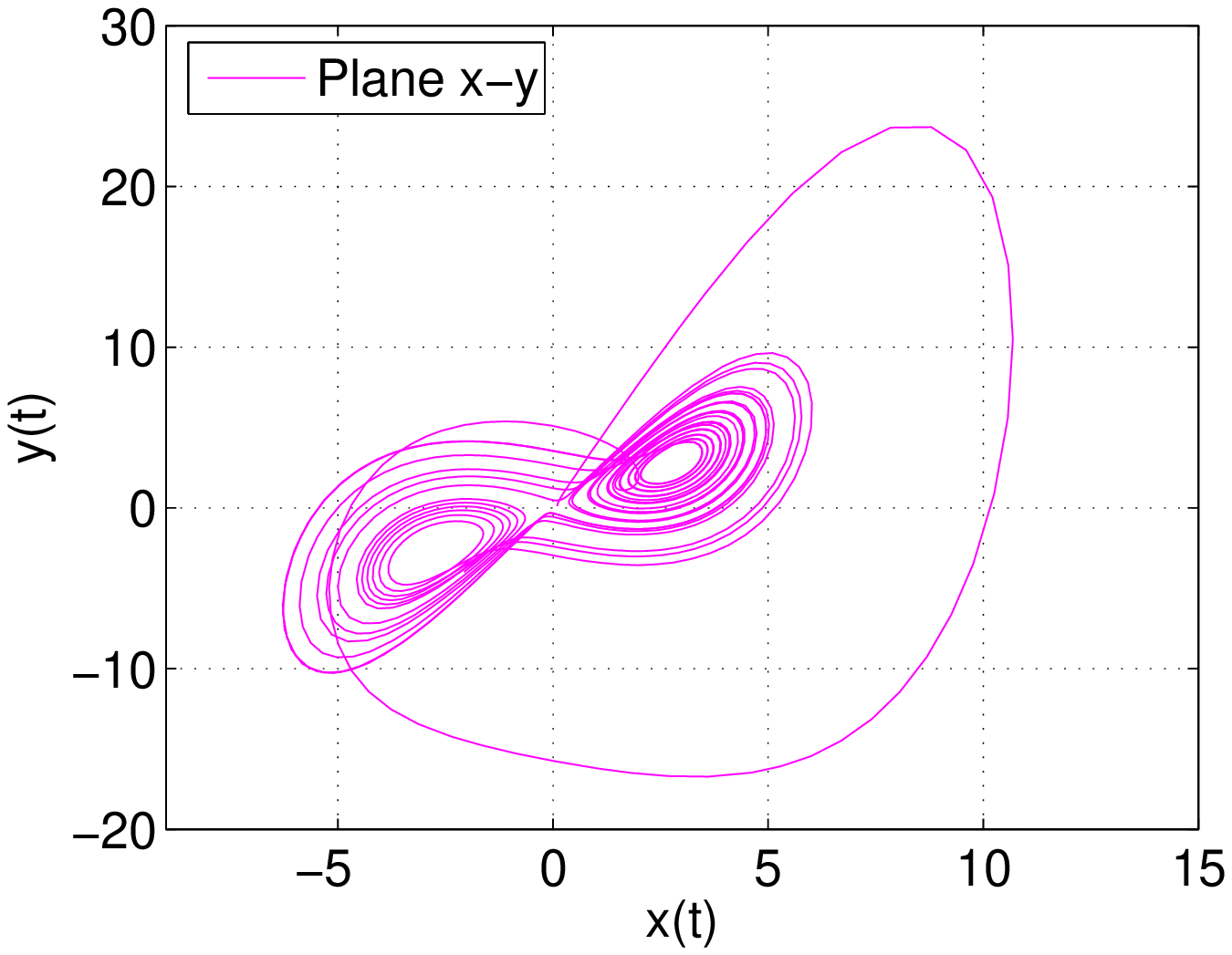}
\caption{The plane x-y with the parameters $a =2$, $D=2/3$
and $\mu =0.9$.}
\label{Fig11}
\end{figure}

The plane x-y of the fractal SL attractor II is given in Figure \ref{Fig11},
where the initial conditions are $x\left( 0 \right)=0.1$, $y\left( 0
\right)=0.1$, and $z\left( 0 \right)=0.1$, and time $t$ changes from $0.1$
to $10^6$.

\begin{figure}[h]
\centering
\includegraphics[width=9cm]{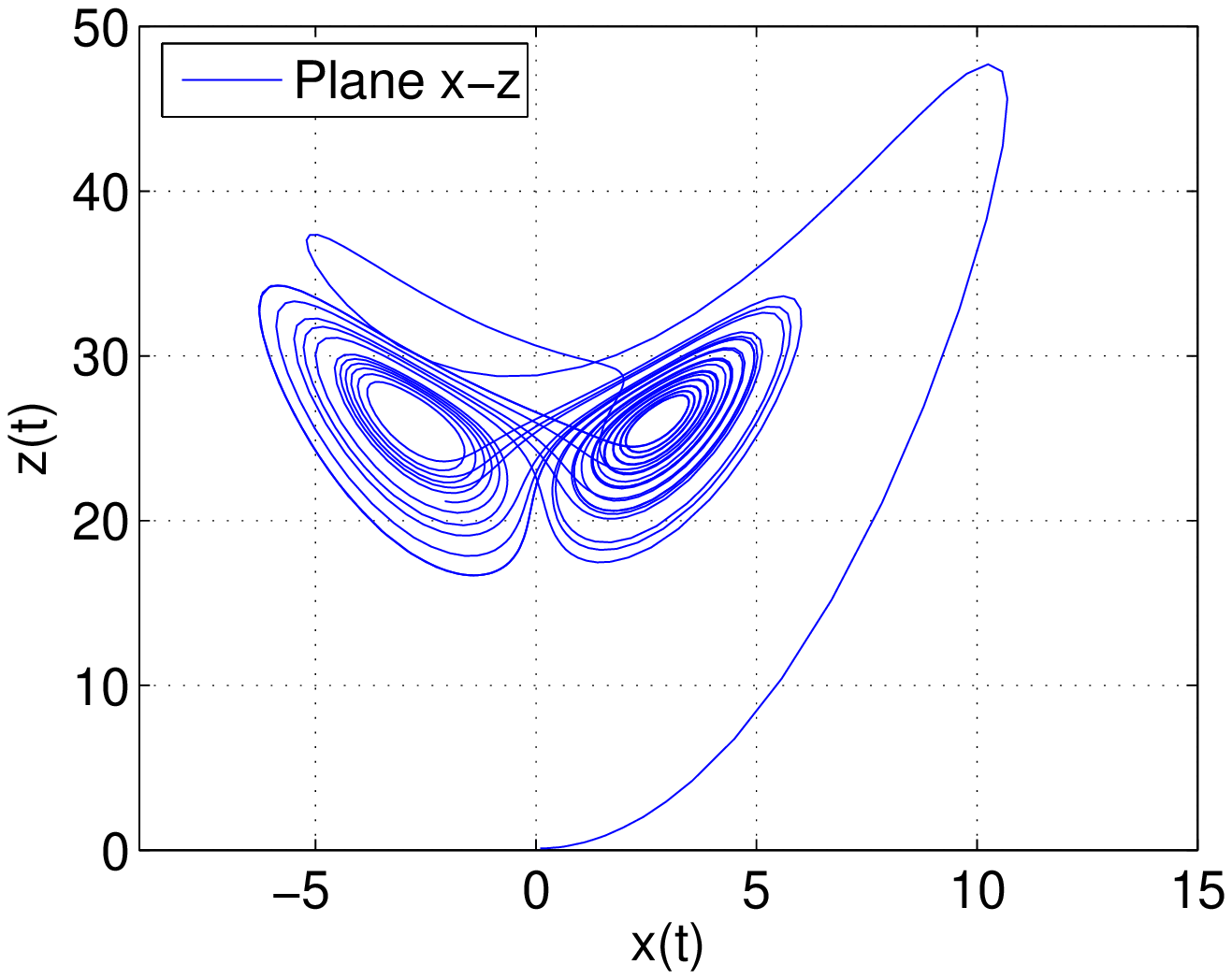}
\caption{The plane x-z with the parameters $a =2$, $D=2/3$
and $\mu =0.9$.}
\label{Fig12}
\end{figure}

The plane x-z of the fractal SL attractor II is given in Figure \ref{Fig12},
where the initial conditions are $x\left( 0 \right)=0.1$, $y\left( 0
\right)=0.1$, and $z\left( 0 \right)=0.1$, and time $t$ changes from $0.1$
to $10^6$.

\begin{figure}[h]
\centering
\includegraphics[width=9cm]{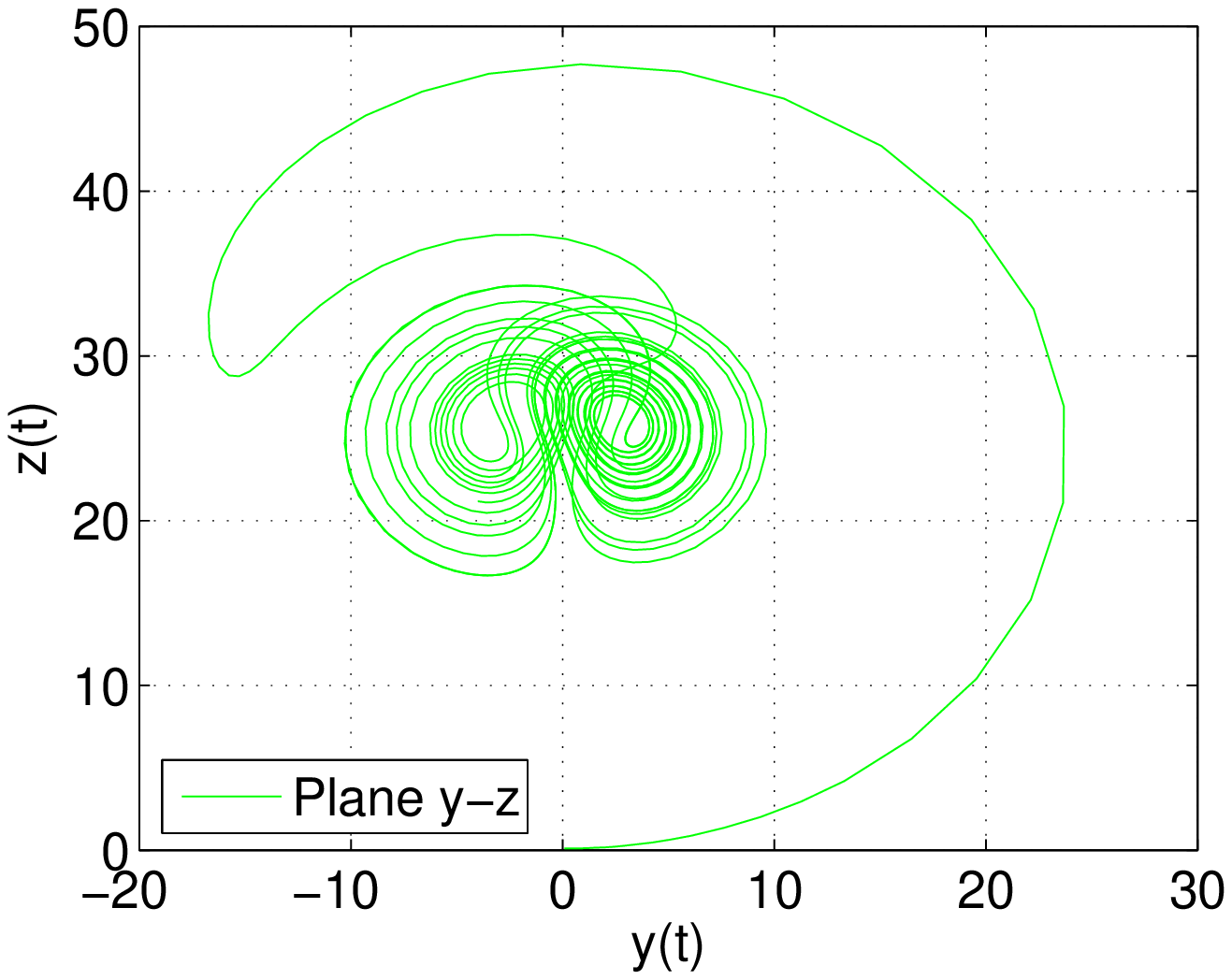}
\caption{The plane y-z with the parameters $a =2$, $D=2/3$
and $\mu =0.9$.}
\label{Fig13}
\end{figure}

The plane y-z of the fractal SL attractor II is depicted in Figure \ref{Fig13},
where the initial conditions are $x\left( 0 \right)=0.1$, $y\left( 0
\right)=0.1$, and $z\left( 0 \right)=0.1$, and time $t$ changes from $0.1$
to $10^6$.

\begin{figure}[h]
\centering
\includegraphics[width=10.5cm]{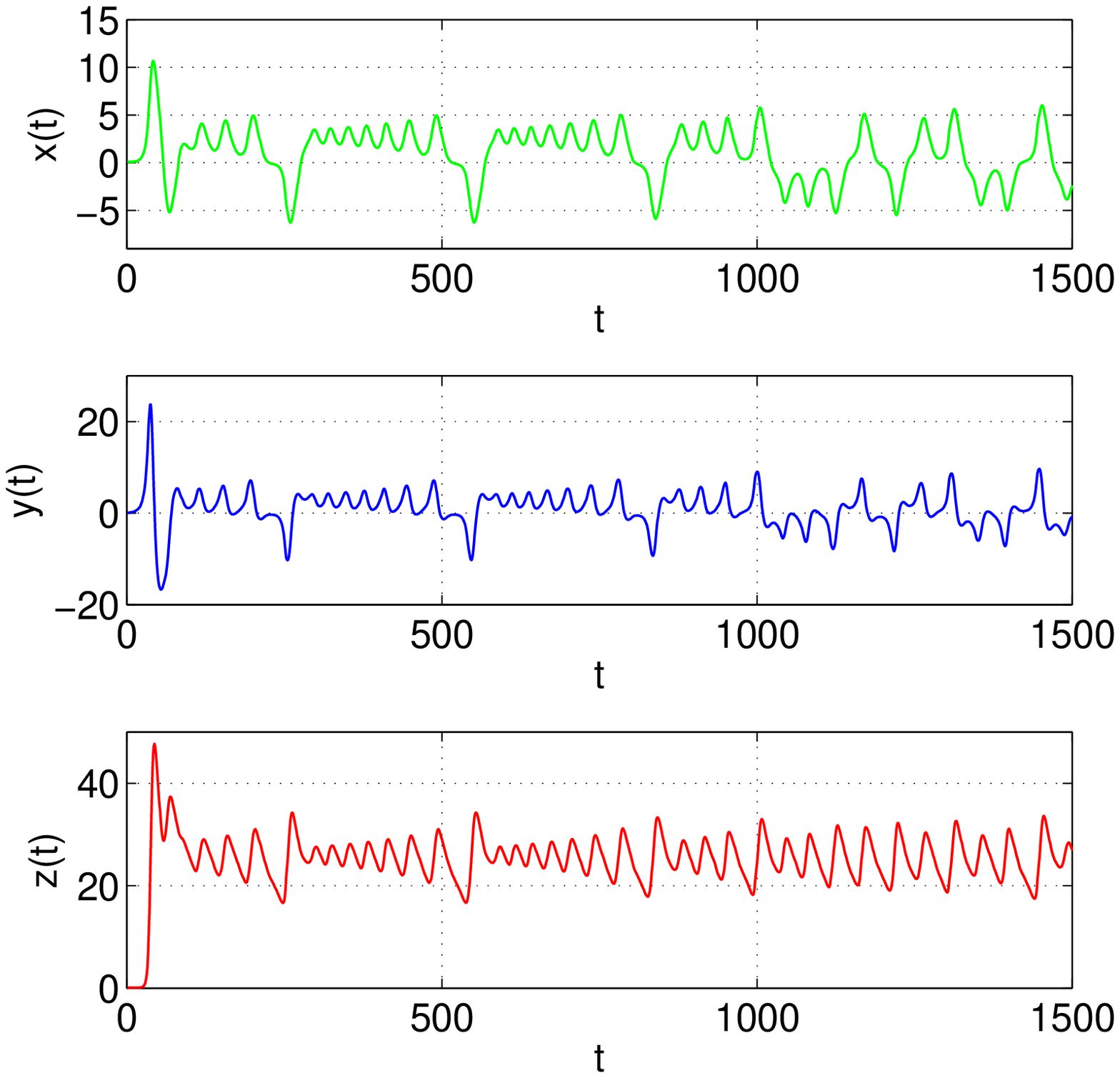}
\caption{The time series for $x\left( t \right)$, $y\left( t \right)$,
and $z\left( t \right)$ with the parameters $a =2$, $D=2/3$
and $\mu =0.9$.}
\label{Fig14}
\end{figure}

The plots of the time series for $x\left( t \right)$, $y\left( t \right)$,
and $z\left( t \right)$ are plotted in Figure \ref{Fig14}, where $D=2/3$, $\mu =0.9$,
and time $t$ changes from $0.1$ to $10^6$.

\section{Comparative results among the Lorenz and fractal SL attractors} \label{Sec:4}

To understand the anomalous behaviors of the fractal SL attractors, we
compare the chaotic trajectories, phase-space parties and time series for
the Lorenz and fractal SL attractors.

\begin{figure}[h]
\centering
\includegraphics[width=9cm]{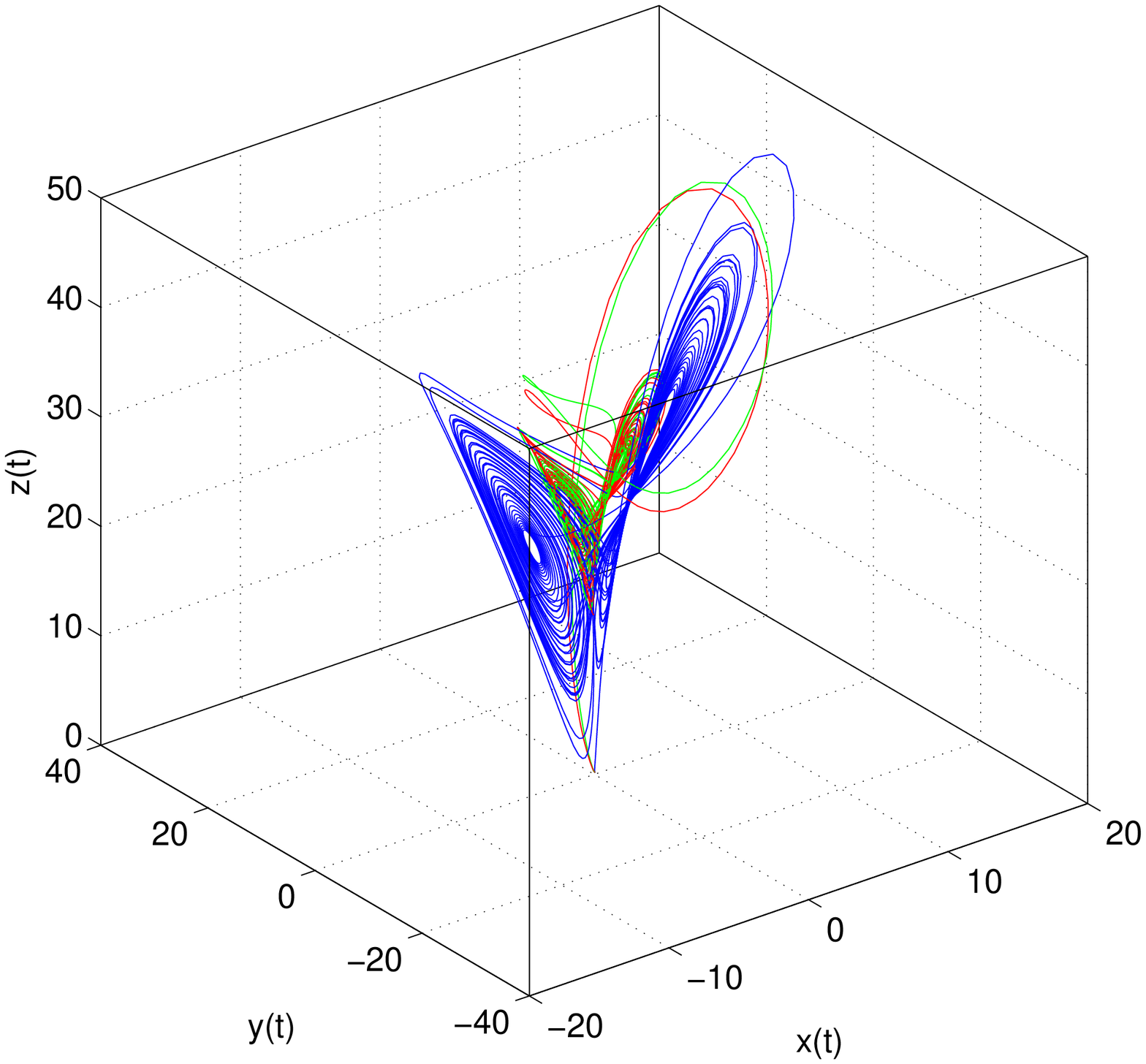}
\caption{The chaotic trajectories for the Lorenz and fractal SL attractors.
The green curve represents the chaotic trajectories for the fractal SL attractor I,
where the initial conditions are $x\left( 0 \right)=0.1$, $y\left( 0
\right)=0.1$, and $z\left( 0 \right)=0.1$, and time $t$ changes from $0.1$
to $10^6$.
The red curve represents the chaotic trajectories for the fractal SL attractor II,
where the initial conditions are $x\left( 0 \right)=0.1$, $y\left( 0
\right)=0.1$, and $z\left( 0 \right)=0.1$, and time $t$ changes from $0.1$
to $10^6$.
The blue curve represents the chaotic trajectories for the Lorenz attractor,
where the initial conditions are $x\left( 0 \right)=0.1$, $y\left( 0
\right)=0.1$, and $z\left( 0 \right)=0.1$, and time $t$ changes from $0$
to $60$.
}
\label{Fig15}
\end{figure}

The chaotic trajectories for the
Lorenz and fractal SL attractors is displayed in Figure \ref{Fig15},
where the initial conditions are $x\left( 0 \right)=0.1$, $y\left( 0
\right)=0.1$, and $z\left( 0 \right)=0.1$.

\begin{figure}[h]
\centering
\includegraphics[width=9cm]{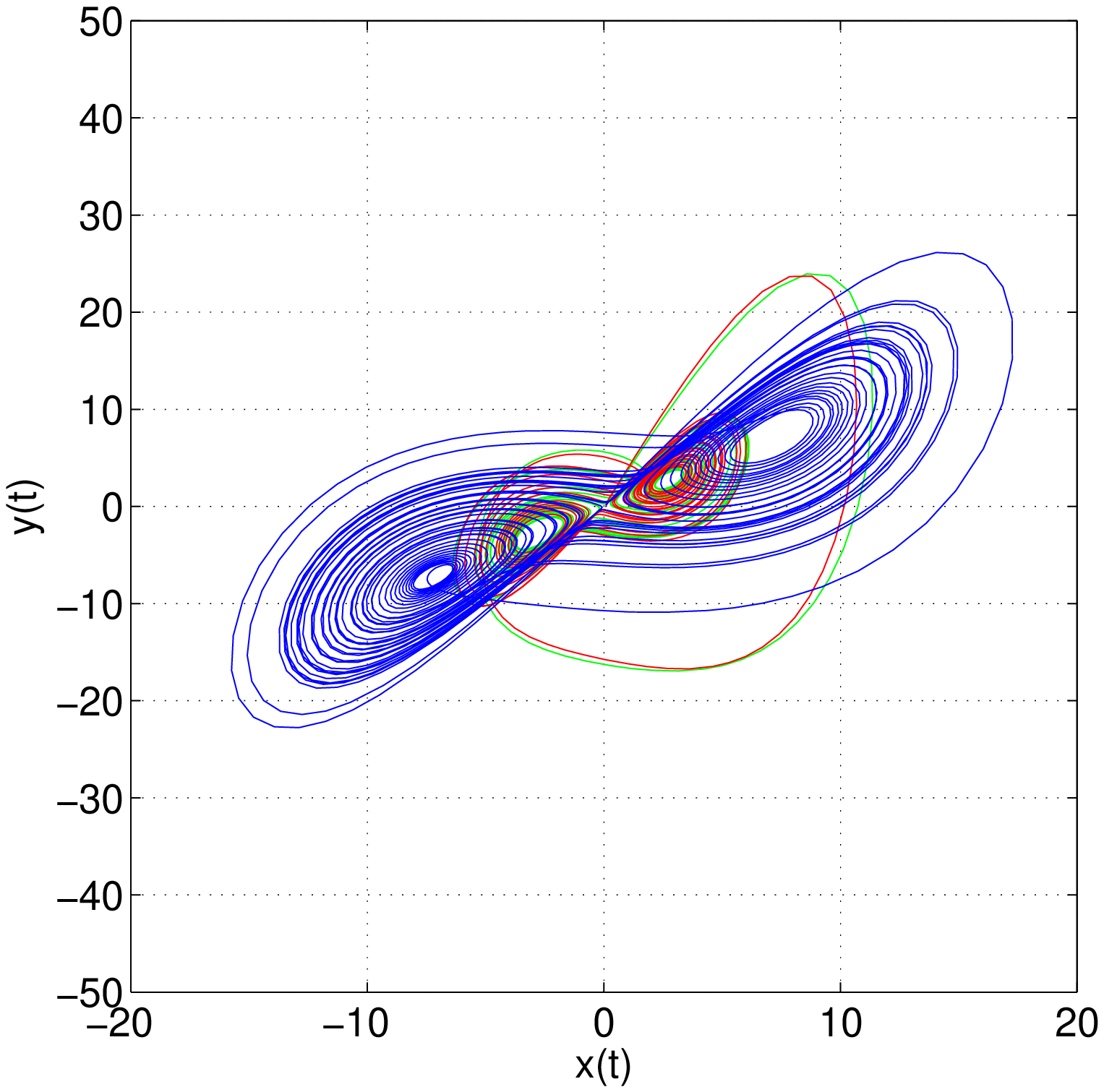}
\caption{The plane x-y for the Lorenz and fractal SL attractors.
The green curve represents the plane x-y for the fractal SL attractor I,
where the initial conditions are $x\left( 0 \right)=0.1$, $y\left( 0
\right)=0.1$, and $z\left( 0 \right)=0.1$, and time $t$ changes from $0.1$
to $10^6$.
The red curve represents the plane x-y for the fractal SL attractor II,
where the initial conditions are $x\left( 0 \right)=0.1$, $y\left( 0
\right)=0.1$, and $z\left( 0 \right)=0.1$, and time $t$ changes from $0.1$
to $10^6$.
The blue curve represents the plane x-y for the Lorenz attractor,
where the initial conditions are $x\left( 0 \right)=0.1$, $y\left( 0
\right)=0.1$, and $z\left( 0 \right)=0.1$, and time $t$ changes from $0$
to $60$.
}
\label{Fig16}
\end{figure}

The plots of the plane x-y for the
Lorenz and fractal SL attractors are demonstrated in Figure \ref{Fig16},
where the initial conditions are $x\left( 0 \right)=0.1$, $y\left( 0
\right)=0.1$, and $z\left( 0 \right)=0.1$.

\begin{figure}[h]
\centering
\includegraphics[width=9cm]{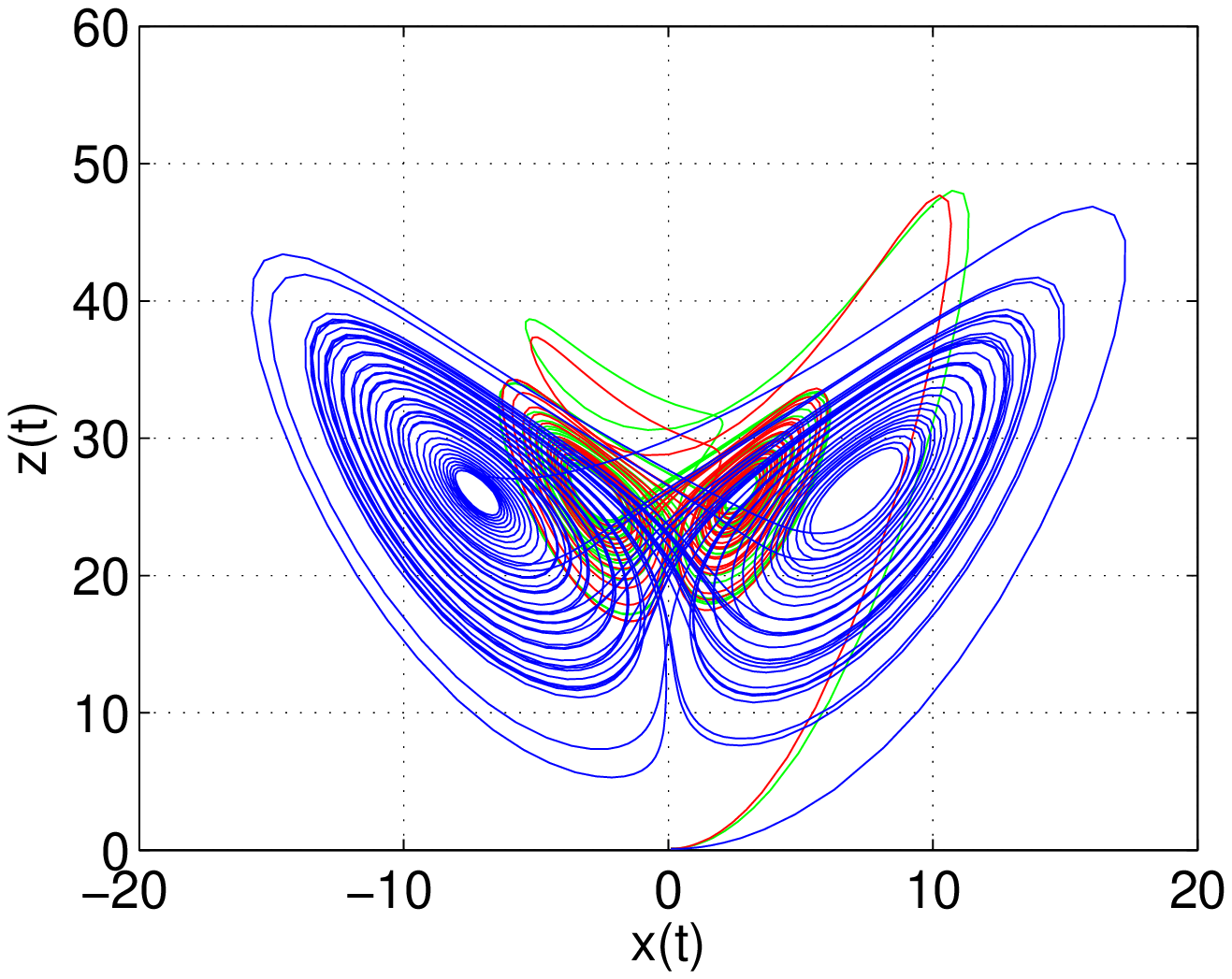}
\caption{The plane x-z for the Lorenz and fractal SL attractors.
The green curve represents the plane x-y for the fractal SL attractor I,
where the initial conditions are $x\left( 0 \right)=0.1$, $y\left( 0
\right)=0.1$, and $z\left( 0 \right)=0.1$, and time $t$ changes from $0.1$
to $10^6$.
The red curve represents the plane x-z for the fractal SL attractor II,
where the initial conditions are $x\left( 0 \right)=0.1$, $y\left( 0
\right)=0.1$, and $z\left( 0 \right)=0.1$, and time $t$ changes from $0.1$
to $10^6$.
The blue curve represents the plane x-z for the Lorenz attractor,
where the initial conditions are $x\left( 0 \right)=0.1$, $y\left( 0
\right)=0.1$, and $z\left( 0 \right)=0.1$, and time $t$ changes from $0$
to $60$.}
\label{Fig17}
\end{figure}

The plots of the plane x-z for the
Lorenz and fractal SL attractors are given in Figure \ref{Fig17},
where the initial conditions are $x\left( 0 \right)=0.1$, $y\left( 0
\right)=0.1$, and $z\left( 0 \right)=0.1$.

The plots of the plane y-z for the
Lorenz and fractal SL attractors are showed in Figure \ref{Fig17},
where the initial conditions are $x\left( 0 \right)=0.1$, $y\left( 0
\right)=0.1$, and $z\left( 0 \right)=0.1$.

\begin{figure}[h]
\centering
\includegraphics[width=9cm]{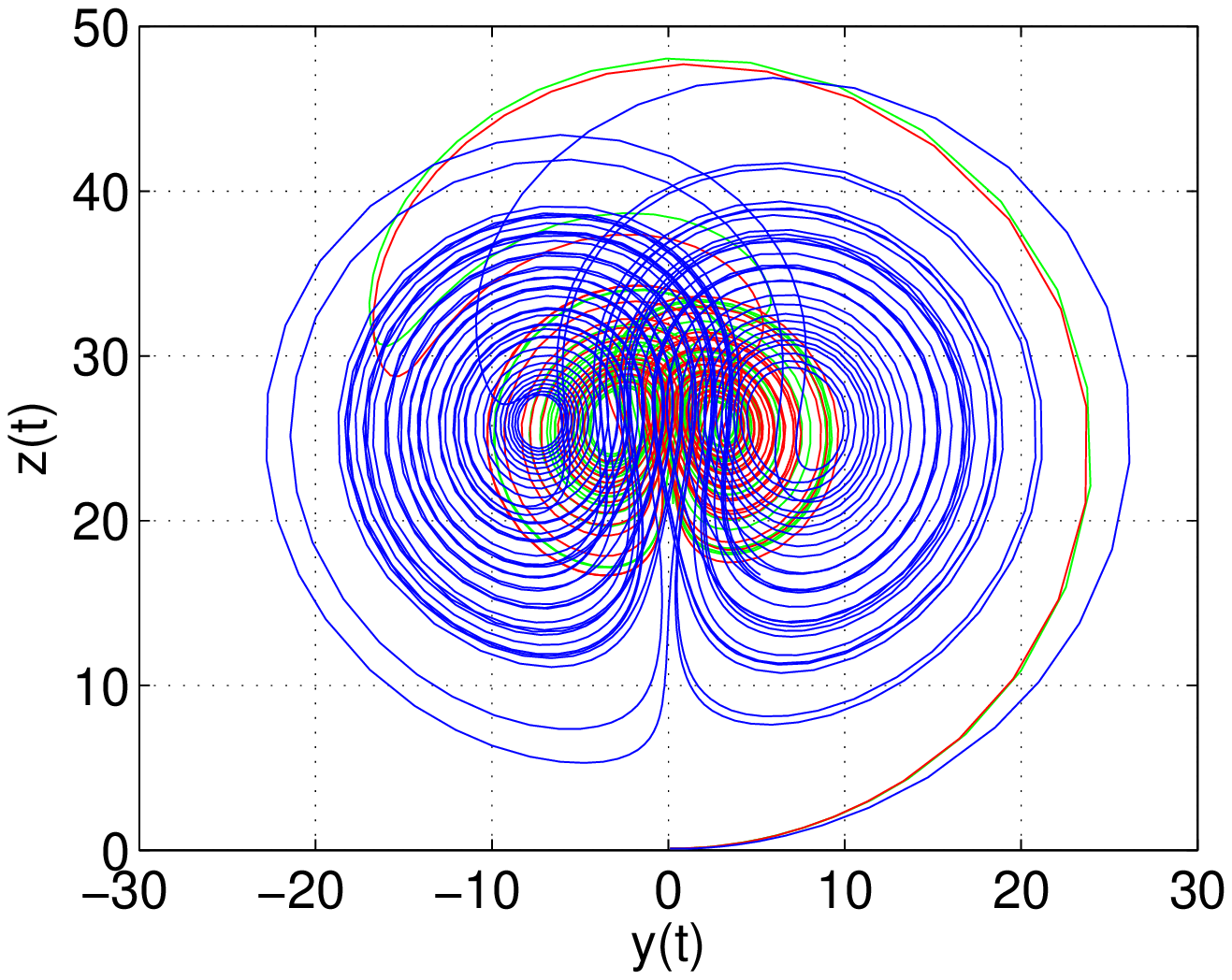}
\caption{The plane y-z for the Lorenz and fractal SL attractors.
The green curve represents the plane y-z for the fractal SL attractor I,
where the initial conditions are $x\left( 0 \right)=0.1$, $y\left( 0
\right)=0.1$, and $z\left( 0 \right)=0.1$, and time $t$ changes from $0.1$
to $10^6$.
The red curve represents the plane y-z for the fractal SL attractor II,
where the initial conditions are $x\left( 0 \right)=0.1$, $y\left( 0
\right)=0.1$, and $z\left( 0 \right)=0.1$, and time $t$ changes from $0.1$
to $10^6$.
The blue curve represents the plane y-z for the Lorenz attractor,
where the initial conditions are $x\left( 0 \right)=0.1$, $y\left( 0
\right)=0.1$, and $z\left( 0 \right)=0.1$, and time $t$ changes from $0$
to $60$.
}
\label{Fig18}
\end{figure}

\begin{figure}[h]
\centering
\includegraphics[width=10.5cm]{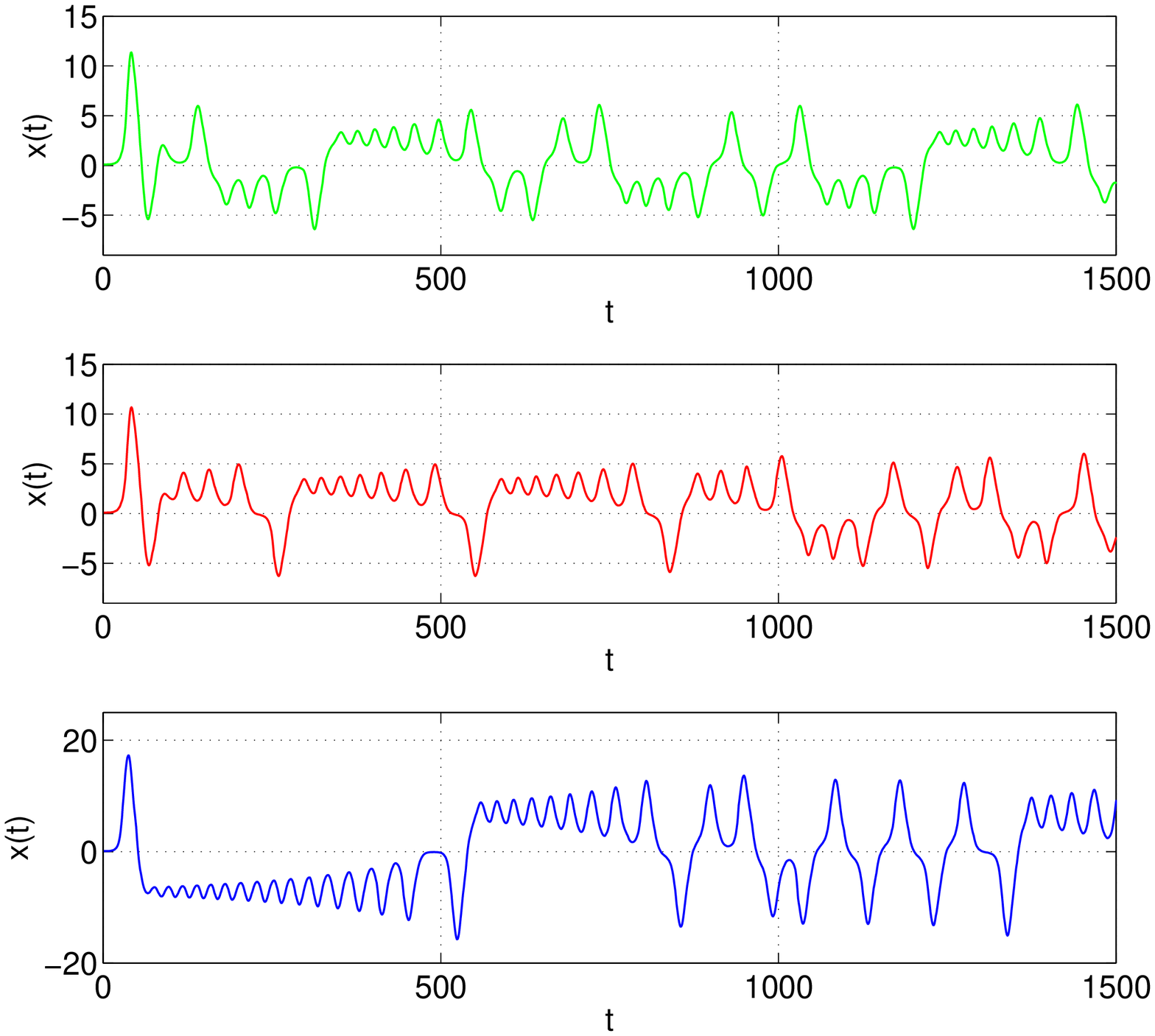}
\caption{The time series $x\left( t \right)$ for the fractal SL and Lorenz attractors.
The green curve represents the time series $x\left( t \right)$ for the fractal SL attractor I,
where the initial conditions are $x\left( 0 \right)=0.1$, $y\left( 0
\right)=0.1$, and $z\left( 0 \right)=0.1$, and time $t$ changes from $0.1$
to $10^6$.
The red curve represents the time series $x\left( t \right)$ for the fractal SL attractor II,
where the initial conditions are $x\left( 0 \right)=0.1$, $y\left( 0
\right)=0.1$, and $z\left( 0 \right)=0.1$, and time $t$ changes from $0.1$
to $10^6$.
The blue curve represents the time series $x\left( t \right)$ for the Lorenz attractor,
where the initial conditions are $x\left( 0 \right)=0.1$, $y\left( 0
\right)=0.1$, and $z\left( 0 \right)=0.1$, and time $t$ changes from $0$
to $60$.
}
\label{Fig19}
\end{figure}

The plots of the time series of $x\left( t \right)$ for
the Lorenz and fractal SL attractors
are plotted in Figure \ref{Fig19}, where the initial conditions are $x\left( 0 \right)=0.1$, $y\left( 0
\right)=0.1$, and $z\left( 0 \right)=0.1$.

\begin{figure}[h]
\centering
\includegraphics[width=10.5cm]{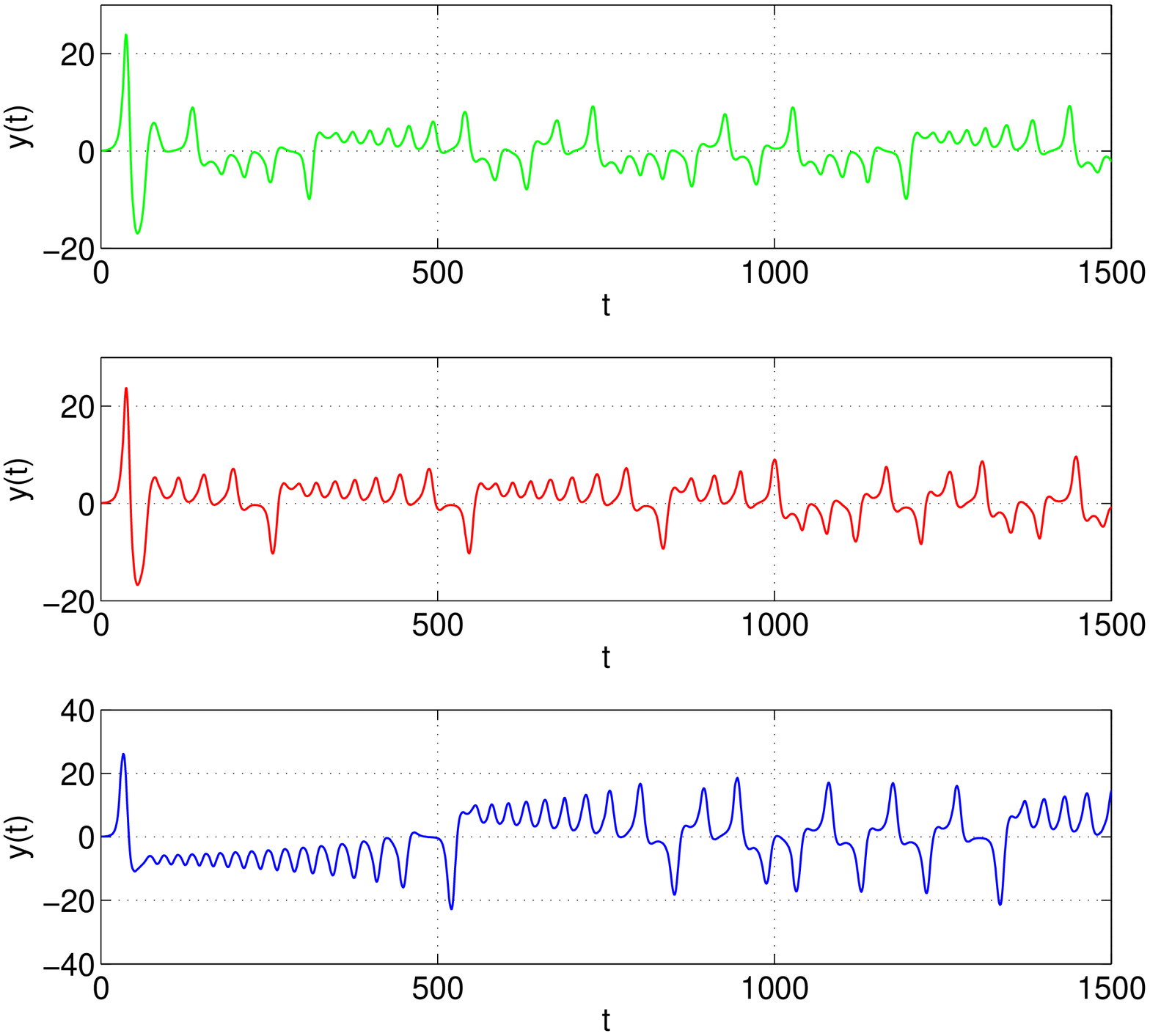}
\caption{The time series $y\left( t \right)$ for the fractal SL and Lorenz attractors.
The green curve represents the time series $y\left( t \right)$ for the fractal SL attractor I,
where the initial conditions are $x\left( 0 \right)=0.1$, $y\left( 0
\right)=0.1$, and $z\left( 0 \right)=0.1$, and time $t$ changes from $0.1$
to $10^6$.
The red curve represents the time series $y\left( t \right)$ for the fractal SL attractor II,
where the initial conditions are $x\left( 0 \right)=0.1$, $y\left( 0
\right)=0.1$, and $z\left( 0 \right)=0.1$, and time $t$ changes from $0.1$
to $10^6$.
The blue curve represents the time series $y\left( t \right)$ for the Lorenz attractor,
where the initial conditions are $x\left( 0 \right)=0.1$, $y\left( 0
\right)=0.1$, and $z\left( 0 \right)=0.1$, and time $t$ changes from $0$
to $60$.
}
\label{Fig20}
\end{figure}

The plots of the time series of $x\left( t \right)$ for
the Lorenz and fractal SL attractors
are plotted in Figure \ref{Fig20}, where the initial conditions are $x\left( 0 \right)=0.1$, $y\left( 0
\right)=0.1$, and $z\left( 0 \right)=0.1$.

\begin{figure}[h]
\centering
\includegraphics[width=10.5cm]{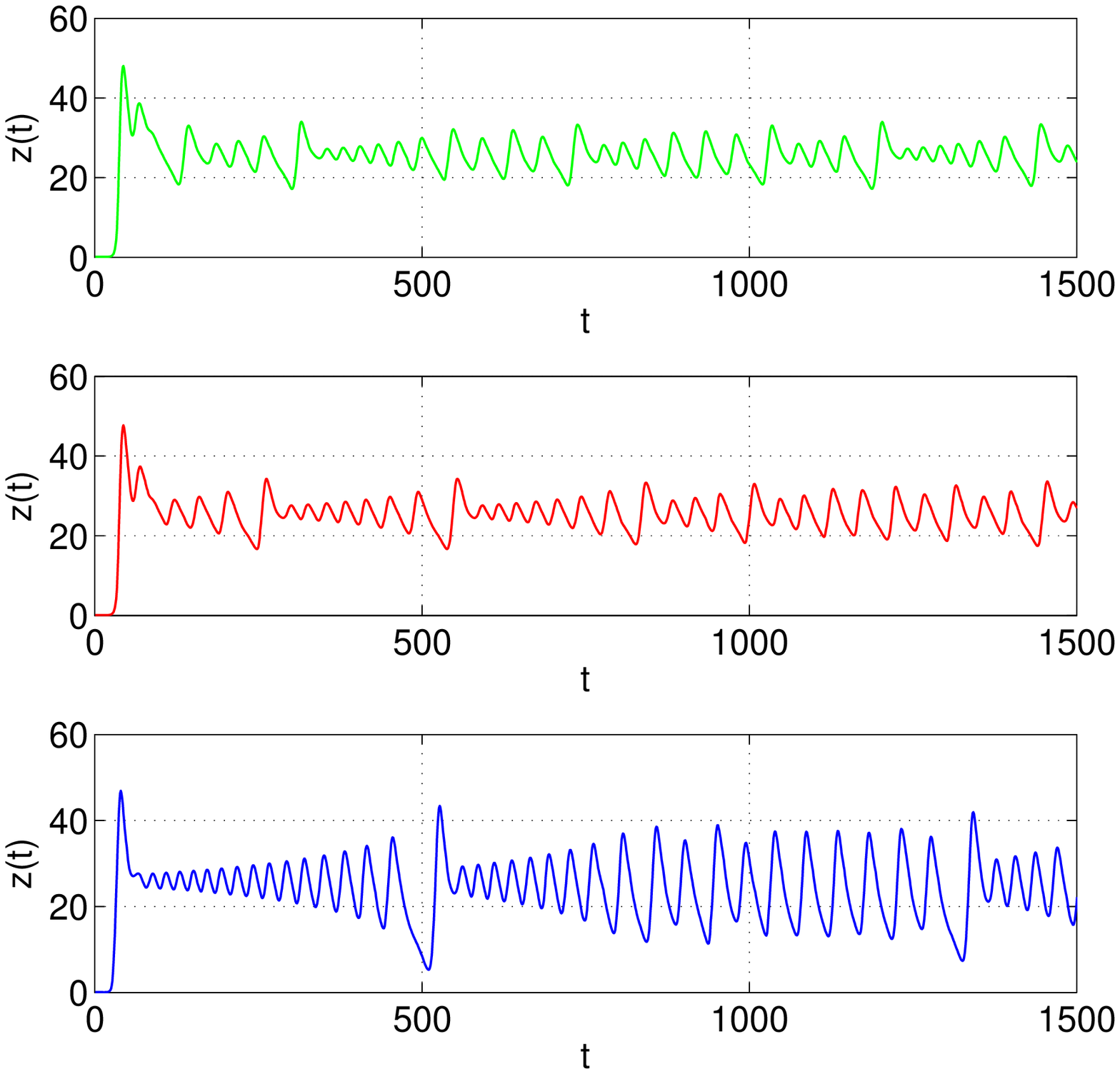}
\caption{The time series $z\left( t \right)$ for the fractal SL and Lorenz attractors.
The green curve represents the time series $y\left( t \right)$ for the fractal SL attractor I,
where the initial conditions are $x\left( 0 \right)=0.1$, $y\left( 0
\right)=0.1$, and $z\left( 0 \right)=0.1$, and time $t$ changes from $0.1$
to $10^6$.
The red curve represents the time series $z\left( t \right)$ for the fractal SL attractor II,
where the initial conditions are $x\left( 0 \right)=0.1$, $y\left( 0
\right)=0.1$, and $z\left( 0 \right)=0.1$, and time $t$ changes from $0.1$
to $10^6$.
The blue curve represents the time series $z\left( t \right)$ for the Lorenz attractor,
where the initial conditions are $x\left( 0 \right)=0.1$, $y\left( 0
\right)=0.1$, and $z\left( 0 \right)=0.1$, and time $t$ changes from $0$
to $60$.
}
\label{Fig21}
\end{figure}

The plots of the time series of $z\left( t \right)$ for
the Lorenz and fractal SL attractors
are plotted in Figure \ref{Fig21}, where the initial conditions are $x\left( 0 \right)=0.1$, $y\left( 0
\right)=0.1$, and $z\left( 0 \right)=0.1$.

It is observed that the plots of the plane x-y for the
fractal SL attractors show the "Wukong effect" (see Figures \ref{Fig6} and \ref{Fig11}) because they look like the face of  "Wukong" who is
from the famous Chinese classical literary work ＾\textit{The Journey to the West}￣ \cite{13}.

Let
\begin{equation}
\label{eq19}
{\rm {\bf X}}=\left( {{\begin{array}{*{20}c}
 {-2} \hfill & 2 \hfill & 0 \hfill \\
 {\frac{3}{10}-z\left( t \right)} \hfill & {-1} \hfill & {-x\left( t
\right)} \hfill \\
 {y\left( t \right)} \hfill & {x\left( t \right)} \hfill & {-27} \hfill \\
\end{array} }} \right)
\end{equation}
and
\begin{equation}
\label{eq20}
{\rm {\bf \Xi }}\left( t \right)=\left( {x\left( t \right),y\left( t
\right),z\left( t \right)} \right).
\end{equation}
The SL system of the SL ordinary differential equations (\ref{eq16}), (\ref{eq17}) and (\ref{eq18})
can be rewritten as
\begin{equation}
\label{eq21}
\frac{10}{3}t^{\frac{2}{3}}\frac{d{\rm {\bf \Xi }}\left( t \right)}{dt}={\rm
{\bf X\Xi }}\left( t \right),
\end{equation}
where
\begin{equation}
\label{eq22}
{\rm {\bf \Xi }}\left( 0 \right)=\left( {x\left( 0 \right),y\left( 0
\right),z\left( 0 \right)} \right)
\end{equation}

Thus, we have the following conjecture:

\textbf{Conjecture}

The chaotic system (\ref{eq21}) with the initial condition (\ref{eq22}) has at lest one
fixed point.

\section{Conclusion and Future Work} \label{Sec:5}
In the present work we have discovered that the SL systems exhibit the chaotic behaviors for
the parameters $a =2$, $D=2/3$, and $\mu =0.9$.
The SL ordinary differential equations were obtained based on the fractal SL derivative
involving the Mandelbrot scaling law. By the comparison among the Lorenz and scaling-law attractors,
it is seen that the fractal SL attractor shows the "Wukong effect" whilst the Lorenz attractor
shows the "Butterfly effect". We proposed the conjecture for the fixed point theory for the fractal SL attractor.
They are nonlinear, non-periodic, three-dimensional and deterministic systems
In the future, we plan to investigate the fractal dimensions
and fixed point theory of the SL attractors. The mathematical structure and applications of the fractal SL system
are also open problems in the study of the SL chaos theory via the fractal SL calculus.

\section*{ACKNOWLEDGMENTS}

This work is supported by the Yue-Qi Scholar of the China University of Mining and Technology (No. 102504180004).

\bibliographystyle{fractals}


\end{document}